\documentclass[11pt,oneside,english]{amsart}
\usepackage{amssymb}
\usepackage{amsmath}
\usepackage{graphicx}
\usepackage{cite}
\usepackage[colorlinks]{hyperref}

\makeatletter \theoremstyle{plain}
\newtheorem{thm}{Theorem}[section]
\newtheorem{lem}[thm]{Lemma}

\newtheorem{corollary}[thm]{Corollary}
\numberwithin{equation}{section}
\numberwithin{figure}{section}
\theoremstyle{definition}
\newtheorem{defn}[thm]{Definition}
\newtheorem{rem}[thm]{Remark}

\newcommand{\bH}{{{\bf H}}}

\newcommand{\C}{{{\mathbb C}}}
\newcommand{\R}{{{\mathbb R}}}

\usepackage{babel}
\makeatother

\begin{document}

\title[Modulus of hyperbolic domains]{Modulus of elementary domains in the hyperbolic plane}

\author{Ioannis D. Platis}
\address{Department of Mathematics, University of Patras, University Campus, GR-26504 Rion, Achaia, Greece.}
\email{idplatis@math.upatras.gr}

\thanks{The author was funded by the MEDICUS programme, No. 83765.}

\begin{abstract}
We study the modulus of curve families inside elementary domains in the hyperbolic plane $\bH^1_\C$. We establish exact expressions for the modulus of connecting and separating curve families in a hyperbolic circular annulus. For the symmetric hyperbolic quadrilateral $Q(a,b)$, we obtain lower bounds by restricting the curve families to particular foliations, and upper bounds by constructing admissible densities. These estimates bracket the modulus. We also present numerical evidence showing how close the bounds are, how they behave asymptotically, and we evaluate them for specific values of $(a,b)$.
\end{abstract}

\maketitle

\section{Introduction}\label{sec-intro}

The conformal modulus of a curve family, equivalently the reciprocal of its extremal length, introduced by Ahlfors and Beurling~\cite{Ahl1, Ahl2}, is a fundamental invariant in geometric function theory and potential theory. The notion has its roots in the length--area method of Gr\"otzsch and was later developed systematically by Beurling, Ahlfors, Teichm\"uller, and others; see \cite{Ahl2, J, O, Vas, Pl1} for comprehensive accounts. Calculating the exact modulus of a domain is a classical problem that usually depends heavily on the geometric symmetries of the underlying region.

Moduli of curve families in the Euclidean plane have been extensively studied, with important applications to quasiconformal mappings, potential theory, and the geometry of Riemann surfaces \cite{Ahl1, Ahl2, Vas, O}. In the hyperbolic setting, the exponential growth of area with radius fundamentally alters the analytic and geometric constraints on curve families, and the theory exhibits features not present in the Euclidean case \cite{Ahl2, Vas, BFP, BBP}. For the hyperbolic annulus and related domains, modulus computations are closely related to the Teichm\"uller theory of Riemann surfaces and the theory of Jenkins--Strebel quadratic differentials \cite{J, Vas, Ahl2, Pl1}. For quadrilaterals, the modulus is a conformal invariant that coincides with the extremal length of the family of curves joining opposite sides, and its determination is equivalent to solving a Dirichlet--Neumann boundary value problem \cite{Ahl2, Vas, N, O, H}. In the case of circular quadrilaterals, exact computations of the modulus typically involve accessory parameters of the associated Schwarzian differential equation, whose determination is a well-known transcendental problem \cite{A-N, N, H}.

In this paper we study the modulus of curve families in two elementary domains of the hyperbolic plane $\bH^1_\C$: the hyperbolic circular annulus $A_R$, see Section~\ref{sec-a}, and the symmetric hyperbolic quadrilateral $Q(a,b)$, see Section~\ref{sec-q}. For the hyperbolic annulus, continuous rotational symmetry simplifies the governing equations into a one-dimensional problem; in this respect our analysis parallels the symmetric criticality principle of Palais \cite{Pa} and the modulus computations for revolution rings in the Heisenberg group \cite{Pl1, BFP, BBP}. This allows us to find exact expressions in closwed form for both the connecting and separating moduli of $A_R$. To the best of the author's knowledge, the explicit formulas of Theorems~\ref{thm-a-j} and~\ref{thm-a-s}, derived here via hyperbolic polar coordinates, have not appeared in this form in the literature; the author has not found them in the references cited herein. While the same moduli can be obtained from the standard conformal coordinate $u=\log\tanh(r/2)$, the self-contained derivation given below is included for completeness and because it illustrates the method that we subsequently apply to the quadrilateral.

For the quadrilateral $Q(a,b)$ the exact modulus is more delicate. We obtain lower bounds by restricting the curve families to particular foliations, and upper bounds by constructing admissible densities. These estimates bracket the modulus. The explicit analytical bounds for $Q(a,b)$ established in Section~\ref{sec-q} appear to be new; to the best of the author's knowledge, they have not appeared in this form in the literature, and the author has not found them in the references cited herein. We do not claim that the lower bounds obtained from these foliations are sharp for the full connecting families.

The paper is organised as follows. Section~\ref{sec-prel} establishes the geometric framework of the Poincar\'e right half-plane model, recalls the definition of curve moduli, and records the reciprocity relation for quadrilateral moduli that we use repeatedly in Section~\ref{sec-q}. Section~\ref{sec-a} treats the hyperbolic circular annulus $A_R$ and derives the exact expressions in closed form for its moduli. Section~\ref{sec-q} introduces the boundary geometry of the quadrilateral $Q(a,b)$ and proves explicit lower and upper bounds. Section~\ref{sec-num} provides numerical evidence comparing the bounds for specific values and, in Theorem~\ref{thm-asymptotic}, records the asymptotic sharpness of the bounds in the limits $b\to0$ and $a\to\infty$. Section~\ref{sec-o} contains concluding remarks.

The author would like to thank Zolt\'an M. Balogh for stimulating discussions.

\section{Preliminaries}\label{sec-prel}

In Section~\ref{sec-hyplane} we state the basic features of the model of the hyperbolic plane we use throughout this paper. In Section~\ref{sec-hpc} we give a detailed description of hyperbolic polar coordinates, see also \cite{Ma}. We mainly use these coordinates in Section~\ref{sec-a}. In Section~\ref{sec-modulusdefn} we define the modulus of curve families inside domains of the hyperbolic plane and explain its equality with the classical modulus of curves in the Euclidean plane. In Section~\ref{sec-recip} we record the reciprocity relation for quadrilateral moduli, which we shall apply in Section~\ref{sec-q}.

\subsection{Hyperbolic plane}\label{sec-hyplane}

We consider the hyperbolic plane $\bH^1_\C$ with coordinates $z=\lambda+it$, $\lambda>0$, $t\in\R$. The metric tensor is given by
\[
g_h=\frac{d\lambda^2+dt^2}{\lambda^2}.
\]
The hyperbolic distance $d_h(z_1,z_2)$ between two points $z_i=\lambda_i+it_i\in\bH^1_\C$, $i=1,2$, is given by
\[
\cosh(d_h(z_1,z_2))=1+\frac{|z_1-z_2|^2}{2\lambda_1\lambda_2}.
\]
The geodesics of $\bH^1_\C$ are either straight lines orthogonal to the boundary line $\lambda=0$ or semicircles centred on that boundary line. Hyperbolic circles $C_h(z_0,R)$ with hyperbolic centre $z_0=\lambda_0+it_0$ and hyperbolic radius $R>0$ are Euclidean circles $C_e((\lambda_0\cosh R,t_0),\lambda_0\sinh R)$.

The line and area elements are, respectively,
\[
ds_h=\frac{\sqrt{d\lambda^2+dt^2}}{\lambda},\qquad d\mathcal{A}_h=\frac{d\lambda\,dt}{\lambda^2}.
\]
The group of orientation-preserving isometries of $(\bH^1_\C,d_h)$ comprises M\"obius transformations of the form
\[
f(z)=\frac{az+ib}{icz+d},\qquad a,b,c,d\in\R,\ ad+bc=1,
\]
that is, it is isomorphic to
\[
\operatorname{PSU}(1,1)=\operatorname{SU}(1,1)/\{\pm I\},\qquad
\operatorname{SU}(1,1)=\left\{\left(\begin{matrix}
a&ib\\
ic&d
\end{matrix}\right),\quad a,b,c,d\in\R,\ ad+bc=1\right\}.
\]

\subsection{Hyperbolic polar coordinates}\label{sec-hpc}

We imitate the standard construction of Euclidean polar coordinates. In the hyperbolic case, the geodesics through $(1,0)$ are either the straight line $t=0$ or semicircles centred on the boundary line $\lambda=0$. We fix the positive real axis and, for $r>0$, consider the point $(e^r,0)$. Every orientation-preserving hyperbolic isometry fixing $1$ has the form
\[
f(z)=\frac{az-ib}{-ibz+a},\qquad a,b\in\R,\ a^2+b^2=1.
\]
Equivalently, we may write
\[
f_\theta(z)=\frac{\cos(\theta/2)z-i\sin(\theta/2)}{-i\sin(\theta/2)z+\cos(\theta/2)},\qquad \theta\in\R.
\]
We now take
\begin{align*}
z=f_\theta(e^r)
&=\frac{e^r\cos(\theta/2)-i\sin(\theta/2)}{-ie^r\sin(\theta/2)+\cos(\theta/2)}\\
&=\frac{2e^r}{(1+e^{2r})+(1-e^{2r})\cos\theta}
+i\frac{\sin\theta(e^{2r}-1)}{(1+e^{2r})+(1-e^{2r})\cos\theta}.
\end{align*}
Therefore,
\[
(\lambda,t)=\left(\frac{1}{\cosh r-\cos \theta\sinh r},\,
\frac{\sin\theta\sinh r}{\cosh r-\cos\theta\sinh r}\right).
\]

\begin{defn}\label{defn-hyppolar}
The map $\Phi:(0,\infty)\times[0,2\pi)\to\bH^1_\C$ given by
\[
(r,\theta)\mapsto\left(\frac{1}{\cosh r-\cos\theta\sinh r},\;
\frac{\sin\theta\sinh r}{\cosh r-\cos\theta\sinh r}\right),
\]
is the \emph{hyperbolic polar coordinates map} in the hyperbolic plane. For $r=0$ we set $\Phi(0,\theta)=(1,0)$ for every $\theta$; the induced map on $[0,\infty)\times[0,2\pi)$ is not injective at $r=0$.
\end{defn}

To find the inverse $\Phi^{-1}$ on $\bH^1_\C\setminus\{(1,0)\}$, we calculate for each $(\lambda,t)\neq(1,0)$:
\[
\frac{\lambda^2+t^2+1}{2\lambda}=\cosh r.
\]
Therefore
\begin{equation}\label{eq-r0}
r=\operatorname{arccosh}\left(\frac{\lambda^2+t^2+1}{2\lambda}\right),
\end{equation}
and the angle $\theta$ is defined uniquely in $[0,2\pi)$ as follows. Put
\[
q=\frac{2t}{\lambda^2+t^2-1}
\]
whenever $\lambda^2+t^2\neq1$, and let $\arctan$ take values in $(-\pi/2,\pi/2)$. Then
\begin{equation}\label{eq-theta0}
\theta=\begin{cases}
0 & \text{if } t=0,\ \lambda\ge1,\\[2mm]
\arctan q & \text{if } \lambda^2+t^2-1>0,\ t>0,\\[2mm]
\pi/2 & \text{if } \lambda^2+t^2=1,\ t>0,\\[2mm]
\pi+\arctan q & \text{if } \lambda^2+t^2-1<0,\ t>0,\\[2mm]
\pi & \text{if } t=0,\ \lambda<1,\\[2mm]
\pi+\arctan q & \text{if } \lambda^2+t^2-1<0,\ t<0,\\[2mm]
3\pi/2 & \text{if } \lambda^2+t^2=1,\ t<0,\\[2mm]
2\pi+\arctan q & \text{if } \lambda^2+t^2-1>0,\ t<0.
\end{cases}
\end{equation}
Therefore, to each $(\lambda,t)\in\bH^1_\C\setminus\{(1,0)\}$ we can assign a unique pair $(r,\theta)$ given by equations (\ref{eq-r0}) and (\ref{eq-theta0}). These equations define the inverse of the restriction of $\Phi$ to $r>0$.

Straightforward computations show that for the Jacobian $J_\Phi$ of $\Phi$ we have
\[
J_\Phi=\lambda^2\sinh r.
\]
The expression for the metric tensor is
\begin{equation}\label{eq-mettens-hyppol}
g_h=dr^2+\sinh^2 r\,d\theta^2.
\end{equation}
In hyperbolic polar coordinates we may thus write the line and area elements respectively as
\[
ds_h=\sqrt{dr^2+\sinh^2 r\,d\theta^2},\qquad d\mathcal{A}_h=\sinh r\,dr\,d\theta.
\]
Finally, in hyperbolic polar coordinates, the circle $C_h((1,0),R)$ is given by the equation $r=R$.

\subsection{Modulus of curve families and the Euclidean connection}\label{sec-modulusdefn}

Let $\Gamma$ be a family of curves in $\bH^1_\C$. Its hyperbolic modulus $\operatorname{Mod}_h(\Gamma)$ is defined via the variational problem
\[
\operatorname{Mod}_h(\Gamma)=\inf_{\rho_h \in \operatorname{Adm}(\Gamma)}\iint_{\bH^1_\C}\rho_h^2\,d\mathcal{A}_h,
\]
where $\operatorname{Adm}(\Gamma)$ denotes the space of Borel measurable densities $\rho_h:\bH^1_\C\to [0,\infty]$ that satisfy the admissibility condition
\[
\int_\gamma\rho_h\,ds_h\ge 1,\qquad \text{for all }\gamma\in\Gamma.
\]
We use the standard convention that densities may be extended-valued; the energy is allowed to be infinite. The energy functional is
\[
\mathcal{E}_h:\operatorname{Adm}(\Gamma)\to [0,\infty],\qquad
\mathcal{E}_h(\rho_h)=\iint_{\bH^1_\C}\rho_h^2\,d\mathcal{A}_h.
\]
We next explain the relationship between the hyperbolic modulus and the Euclidean modulus $\operatorname{Mod}_e(\Gamma)$. Because the hyperbolic metric
\[
g_h=\lambda^{-2}(d\lambda^2+dt^2)
\]
is conformally equivalent to the standard Euclidean metric $g_e=d\lambda^2+dt^2$, the conformal invariance of the modulus implies that for any curve family $\Gamma$,
\[
\operatorname{Mod}_h(\Gamma)=\operatorname{Mod}_e(\Gamma).
\]
This identity can be seen directly by mapping any admissible hyperbolic density $\rho_h$ to a corresponding Euclidean density via $\rho_e(\lambda,t)=\rho_h(\lambda,t)/\lambda$. Under this transformation, both the admissibility line integrals and the energy area integrals are equal:
\[
\int_\gamma \rho_h\,ds_h
=\int_\gamma \left(\frac{\rho_h}{\lambda}\right)\sqrt{d\lambda^2+dt^2}
=\int_\gamma \rho_e\,ds_e,
\]
\[
\iint \rho_h^2\,d\mathcal{A}_h
=\iint \left(\frac{\rho_h}{\lambda}\right)^2 d\lambda\,dt
=\iint \rho_e^2\,d\mathcal{A}_e.
\]
Consequently, computing $\operatorname{Mod}_h(\Gamma)$ is equivalent to computing the Euclidean modulus of the same curve family.

\subsection{The reciprocity relation for quadrilateral moduli}\label{sec-recip}

We record, with proof, the reciprocity relation that we invoke in Section~\ref{sec-q}. Although the relation is classical, we include the short argument based on conformal invariance so that the paper is self-contained.

\begin{lem}[Reciprocity]\label{lem-recip}
Let $Q$ be a topological quadrilateral with distinguished boundary arcs $S_1,S_2,S_3,S_4$ in cyclic order. Let $\Gamma_1$ be the family of curves in $Q$ joining $S_1$ to $S_3$, and let $\Gamma_2$ be the family joining $S_2$ to $S_4$. Then
\[
\operatorname{Mod}_h(\Gamma_1)\operatorname{Mod}_h(\Gamma_2)=1.
\]
\end{lem}

\begin{proof}
By the classical conformal-mapping theorem for quadrilaterals (see Chapter~I of \cite{Ahl2} and Chapter~VI of \cite{N}), there exists a conformal map $f$ from $Q$ onto a Euclidean rectangle
\[
R=\{(x,y):0<x<a,\ 0<y<b\},
\]
which sends the four vertices of $Q$ to the four corners of $R$, the arcs $S_1$ and $S_3$ to the horizontal sides $y=0$ and $y=b$, and the arcs $S_2$ and $S_4$ to the vertical sides $x=0$ and $x=a$. Because the modulus is conformally invariant, it suffices to verify the relation for the Euclidean rectangle $R$.

On $R$ the family $\Gamma_1$ corresponds to the horizontal family joining $y=0$ to $y=b$, and $\Gamma_2$ to the vertical family joining $x=0$ to $x=a$. For the horizontal family, the constant density $\rho_1=1/b$ is admissible, since any horizontal curve has Euclidean length at least $b$ in the vertical direction; its energy is
\[
\mathcal{E}_e(\rho_1)=\iint_R \frac{1}{b^2}\,dx\,dy=\frac{ab}{b^2}=\frac{a}{b}.
\]
For the vertical family, the constant density $\rho_2=1/a$ is admissible, and its energy is
\[
\mathcal{E}_e(\rho_2)=\iint_R \frac{1}{a^2}\,dx\,dy=\frac{ab}{a^2}=\frac{b}{a}.
\]
Both densities are extremal: the associated foliations are the horizontal and vertical segments, which are orthogonal, and the standard length--area argument (or the equivalent Dirichlet-energy argument) shows that the infimum of the energy is attained by them. Hence
\[
\operatorname{Mod}_e(\Gamma_1)=\frac{a}{b},\qquad
\operatorname{Mod}_e(\Gamma_2)=\frac{b}{a}.
\]
By the conformal invariance $\operatorname{Mod}_h=\operatorname{Mod}_e$ established in Section~\ref{sec-modulusdefn}, the same equalities hold for the hyperbolic moduli. Multiplying the two equalities yields $\operatorname{Mod}_h(\Gamma_1)\operatorname{Mod}_h(\Gamma_2)=1$.
\end{proof}

\begin{rem}
In the notation of Section~\ref{sec-q}, the two connecting families $\Gamma_Q^1$ and $\Gamma_Q^2$ of the quadrilateral $Q(a,b)$ are related by Lemma~\ref{lem-recip}:
\[
\operatorname{Mod}_h(\Gamma_Q^1)\operatorname{Mod}_h(\Gamma_Q^2)=1.
\]
This is the form in which we apply the reciprocity relation.
\end{rem}

\section{Modulus of the hyperbolic annulus}\label{sec-a}

\begin{figure}[!h]
\centering
\includegraphics[scale=0.6]{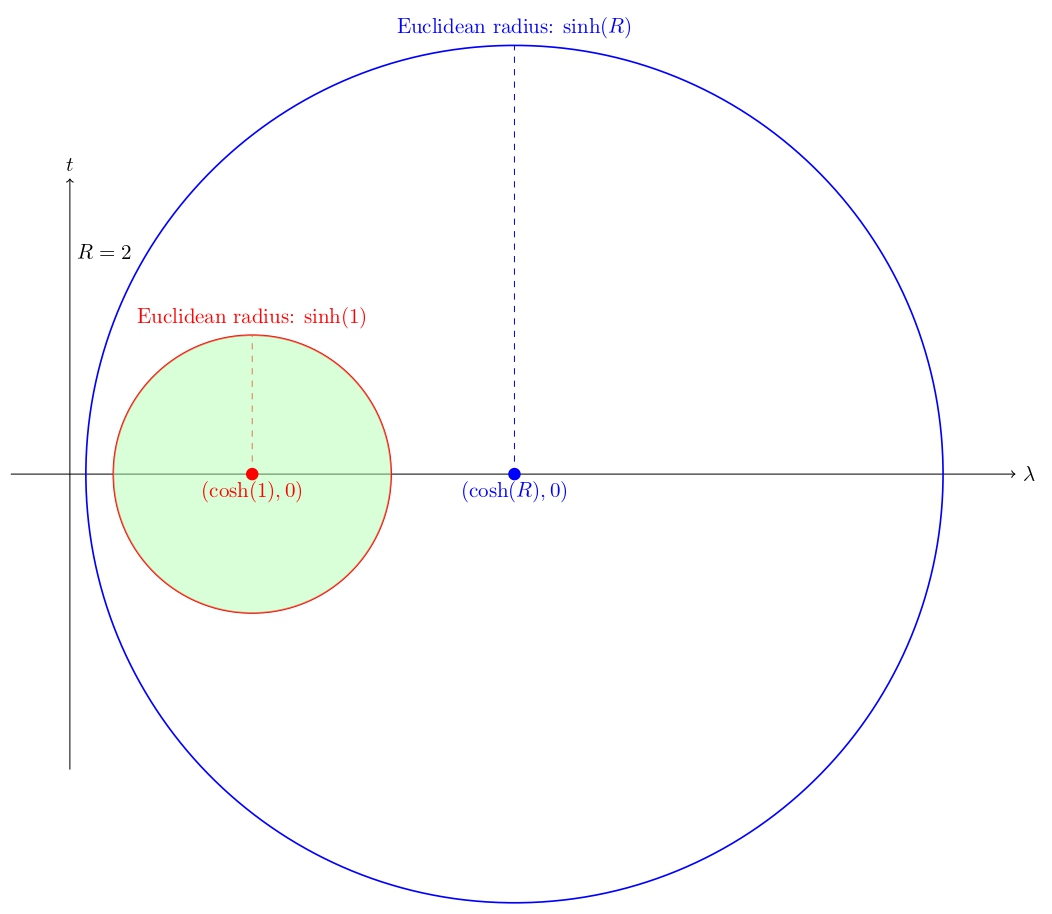}
\caption{The white region between the circles is the hyperbolic annulus.}
\label{fig:3.1}
\end{figure}

The hyperbolic annulus, see Figure~\ref{fig:3.1}, possesses perfect rotational symmetry; the underlying potential-theoretic problems reduce to solving one-dimensional ordinary differential equations. This structure allows us to determine the extremal densities and state the main explicit formulae of the paper.

\begin{defn}
For $R>1$, the hyperbolic circular annulus $A_R$ is defined in hyperbolic polar coordinates as
\[
A_R=\{\Phi(r,\theta): 1<r<R,\ 0\le\theta<2\pi\}.
\]
\end{defn}

The following theorem gives the modulus for the family of curves connecting the two boundary components of the annulus.

\begin{thm}\label{thm-a-j}
Let $\Gamma_A^1$ be the family of rectifiable curves $\gamma:[0,1]\to\overline{A_R}$ such that $\gamma((0,1))\subset A_R$ and the endpoints lie on different boundary components of $A_R$. Then
\[
\operatorname{Mod}_h(\Gamma_A^1)=\frac{2\pi}{\log\left(\frac{\tanh(R/2)}{\tanh(1/2)}\right)}.
\]
\end{thm}

\begin{proof}
Consider the radial subfamily $\Gamma_A'\subset\Gamma_A^1$ composed of the segments
\[
\gamma_\theta(r)=\Phi(r,\theta),\qquad 1\le r\le R,\quad \theta\in[0,2\pi).
\]
Every $\rho\in\operatorname{Adm}(\Gamma_A^1)$ is admissible for the radial subfamily $\Gamma_A'$. Hence, for each $\theta\in[0,2\pi)$,
\[
\int_{1}^{R}\rho(r,\theta)\,dr\ge1.
\]
Applying the Cauchy--Schwarz inequality gives
\[
1\le\left(\int_{1}^{R}\rho(r,\theta)\,dr\right)^2
\le\left(\int_{1}^{R}\rho^2(r,\theta)\sinh r\,dr\right)
\left(\int_{1}^{R}\frac{dr}{\sinh r}\right).
\]
Integrating over $\theta\in[0,2\pi)$ yields
\[
\iint_{A_R}\rho^2\,d\mathcal{A}_h
\ge 2\pi\left[\log\left(\frac{\tanh(R/2)}{\tanh(1/2)}\right)\right]^{-1}.
\]
Taking the infimum over all admissible densities proves that this value is a lower bound for $\operatorname{Mod}_h(\Gamma_A^1)$.

To prove that it is also an upper bound, define the density
\[
\rho_{A,1}(r,\theta)=\frac{c_0}{\sinh r}\chi_{A_R}(r,\theta),
\qquad
c_0=\left[\log\left(\frac{\tanh(R/2)}{\tanh(1/2)}\right)\right]^{-1}.
\]
For any curve $\gamma\in\Gamma_A^1$ connecting $r=1$ to $r=R$, oriented so that $r(0)=1$ and $r(1)=R$, we have $ds_h\ge |dr|$. Thus
\[
\int_\gamma \rho_{A,1}\,ds_h
\ge c_0\int_\gamma \frac{|dr|}{\sinh r}
\ge c_0\left|\int_0^1\frac{r'(s)}{\sinh r(s)}\,ds\right|
=c_0\left|\log\tanh\frac{R}{2}-\log\tanh\frac{1}{2}\right|=1.
\]
Therefore $\rho_{A,1}\in\operatorname{Adm}(\Gamma_A^1)$. Moreover,
\[
\mathcal{E}_h(\rho_{A,1})
=2\pi c_0^2\int_1^R\frac{dr}{\sinh r}
=2\pi c_0.
\]
Hence
\[
\operatorname{Mod}_h(\Gamma_A^1)\le 2\pi c_0.
\]
Combining the lower and upper bounds completes the proof.
\end{proof}

The next theorem addresses the family of closed curves that topologically separate the inner boundary from the outer boundary.

\begin{thm}\label{thm-a-s}
Let $\Gamma_A^2$ be the family of rectifiable closed curves contained in $A_R$ whose images separate the two boundary components of $A_R$. Then
\[
\operatorname{Mod}_h(\Gamma_A^2)=\frac{1}{2\pi}\log\left(\frac{\tanh(R/2)}{\tanh(1/2)}\right).
\]
\end{thm}

\begin{proof}
Let $\Gamma_A''\subset\Gamma_A^2$ be the subfamily of concentric hyperbolic circles
\[
\gamma_r(s)=\Phi(r,s),\qquad 0\le s<2\pi,\quad \gamma_r(2\pi)=\gamma_r(0),
\]
with $r\in(1,R)$. Every $\rho\in\operatorname{Adm}(\Gamma_A^2)$ is admissible for $\Gamma_A''$. Hence
\[
\int_{0}^{2\pi}\rho(r,s)\sinh r\,ds\ge1,
\]
which implies
\[
\frac{1}{\sinh r}\le\int_{0}^{2\pi}\rho(r,s)\,ds.
\]
We integrate this inequality over $r\in(1,R)$ and then apply the Cauchy--Schwarz inequality:
\begin{align*}
\int_{1}^{R}\frac{1}{\sinh r}\,dr
&\le\iint_{A_R}\rho(r,s)(\sinh r)^{1/2}(\sinh r)^{-1/2}\,ds\,dr\\
&\le\left(\iint_{A_R}\rho^2\,d\mathcal{A}_h\right)^{1/2}
\left(2\pi\int_{1}^{R}\frac{1}{\sinh r}\,dr\right)^{1/2}.
\end{align*}
Squaring both sides and isolating the double integral gives
\[
\iint_{A_R}\rho^2\,d\mathcal{A}_h
\ge\frac{1}{2\pi}\int_{1}^{R}\frac{1}{\sinh r}\,dr
=\frac{1}{2\pi}\log\left(\frac{\tanh(R/2)}{\tanh(1/2)}\right).
\]
Thus this value is a lower bound for $\operatorname{Mod}_h(\Gamma_A^2)$.

To obtain an upper bound as well, consider the density
\[
\rho_{A,2}(r,\theta)=\frac{1}{2\pi\sinh r}\chi_{A_R}(r,\theta).
\]
Let $\gamma$ be any closed curve in $A_R$ that separates its boundary components. Then
\[
\int_\gamma |d\theta|\ge 2\pi.
\]
Using $ds_h=\sqrt{dr^2+\sinh^2 r\,d\theta^2}$, we have $ds_h\ge \sinh r\,|d\theta|$. Hence
\[
\int_\gamma \rho_{A,2}\,ds_h
=\int_\gamma \frac{1}{2\pi\sinh r}\,ds_h
\ge \frac{1}{2\pi}\int_\gamma |d\theta|
\ge1.
\]
Therefore $\rho_{A,2}\in\operatorname{Adm}(\Gamma_A^2)$. Finally,
\[
\mathcal{E}_h(\rho_{A,2})
=\frac{1}{2\pi}\int_1^R\frac{dr}{\sinh r}
=\frac{1}{2\pi}\log\left(\frac{\tanh(R/2)}{\tanh(1/2)}\right).
\]
Thus
\[
\operatorname{Mod}_h(\Gamma_A^2)\le
\frac{1}{2\pi}\log\left(\frac{\tanh(R/2)}{\tanh(1/2)}\right),
\]
which completes the proof.
\end{proof}

\section{The hyperbolic quadrilateral and variational bounds}\label{sec-q}

\begin{figure}[!h]
\centering
\includegraphics[scale=0.6]{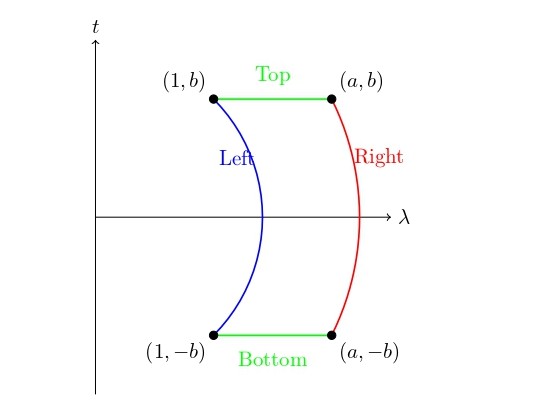}
\caption{A symmetric hyperbolic quadrilateral.}
\label{fig:2.1}
\end{figure}

We consider geodesic quadrilaterals bounded by two horizontal segments and two arcs of Euclidean circles centred at the origin. These quadrilaterals are symmetric with respect to the real axis; see Figure~\ref{fig:2.1}.

\begin{defn}\label{defn-q}
For $a>1$ and $b>0$, define
\[
Q(a,b)=\left\{(\lambda,t)\in\bH^1_\C:
-b<t<b,\quad
\sqrt{1+b^2-t^2}<\lambda<\sqrt{a^2+b^2-t^2}
\right\}.
\]
\end{defn}

The two horizontal sides $t=\pm b$ have equal hyperbolic length $\log a$. The two circular sides have hyperbolic lengths $2\operatorname{arsinh} b$ and $2\operatorname{arsinh}(b/a)$; they are unequal whenever $a>1$. We study the general family $Q(a,b)$ with $a>1$ and $b>0$. For later use we note that the Euclidean area $\mathcal{A}_e(Q(a,b))$ is given by
\begin{equation}\label{area-Q}
\mathcal{A}_e(Q(a,b))
=b(a-1)+(a^2+b^2)\arctan(b/a)-(1+b^2)\arctan b.
\end{equation}

\subsection{Modulus estimate for curves connecting the circular arcs}\label{sec-q-1}

Let $\Gamma_Q^1$ be the family of rectifiable curves inside $Q(a,b)$ connecting the left circular arc $L$ and the right circular arc $R$. We obtain a lower bound by restricting the curve family to the horizontal foliation $\Gamma_Q'\subset\Gamma_Q^1$, whose elements are the horizontal segments
\[
\gamma_t(\lambda)=(\lambda,t),\qquad
\lambda_1(t)=\sqrt{1+b^2-t^2}<\lambda<\lambda_2(t)=\sqrt{a^2+b^2-t^2},
\]
with $t\in(-b,b)$.

\begin{thm}\label{thm-q-c}
Let $a>1$ and $b>0$. Then
\[
\operatorname{Mod}_h(\Gamma_Q^1)\ge \operatorname{Mod}_h(\Gamma_Q')
=\frac{b(a+1)+(a^2+b^2)\arctan(b/a)+(1+b^2)\arctan b}{a^2-1}.
\]
\end{thm}

\begin{proof}
Let $\rho\in\operatorname{Adm}(\Gamma_Q^1)$ be admissible for the full family. Since $\Gamma_Q'\subset\Gamma_Q^1$, the density $\rho$ is admissible for the horizontal foliation. Hence, for every $t\in(-b,b)$,
\[
\int_{\gamma_t}\rho\,ds_h
=\int_{\lambda_1(t)}^{\lambda_2(t)}\rho(\lambda,t)\frac{d\lambda}{\lambda}
\ge1.
\]
By the Cauchy--Schwarz inequality,
\[
1\le\left(\int_{\lambda_1(t)}^{\lambda_2(t)}
\frac{\rho(\lambda,t)}{\lambda}\,d\lambda\right)^2
\le\left(\int_{\lambda_1(t)}^{\lambda_2(t)}
\frac{\rho^2(\lambda,t)}{\lambda^2}\,d\lambda\right)
(\lambda_2(t)-\lambda_1(t)).
\]
Thus, for each $t\in(-b,b)$,
\[
\int_{\lambda_1(t)}^{\lambda_2(t)}
\rho^2(\lambda,t)\frac{d\lambda}{\lambda^2}
\ge\frac{1}{\lambda_2(t)-\lambda_1(t)}.
\]
Integrating over $t\in(-b,b)$ yields
\[
\iint_{Q(a,b)}\rho^2\,d\mathcal{A}_h
\ge\int_{-b}^{b}\frac{dt}{\lambda_2(t)-\lambda_1(t)}.
\]
Taking the infimum over all $\rho\in\operatorname{Adm}(\Gamma_Q^1)$ gives
\[
\operatorname{Mod}_h(\Gamma_Q^1)
\ge\int_{-b}^{b}\frac{dt}{\sqrt{a^2+b^2-t^2}-\sqrt{1+b^2-t^2}}.
\]
Now
\[
\frac{1}{\sqrt{a^2+b^2-t^2}-\sqrt{1+b^2-t^2}}
=\frac{\sqrt{a^2+b^2-t^2}+\sqrt{1+b^2-t^2}}{a^2-1}.
\]
Therefore
\begin{align*}
\int_{-b}^{b}\frac{dt}{\lambda_2(t)-\lambda_1(t)}
&=\frac{1}{a^2-1}
\left(
\int_{-b}^{b}\sqrt{a^2+b^2-t^2}\,dt
+\int_{-b}^{b}\sqrt{1+b^2-t^2}\,dt
\right)\\
&=\frac{b(a+1)+(a^2+b^2)\arctan(b/a)+(1+b^2)\arctan b}{a^2-1}.
\end{align*}
The density
\[
\rho_{\ell,1}(\lambda,t)
=\frac{\lambda}{\lambda_2(t)-\lambda_1(t)}\chi_{Q(a,b)}
\]
is extremal for $\Gamma_Q'$. Indeed,
\[
\int_{\gamma_t}\rho_{\ell,1}\,ds_h=1,
\qquad
\mathcal{E}_h(\rho_{\ell,1})
=\frac{b(a+1)+(a^2+b^2)\arctan(b/a)+(1+b^2)\arctan b}{a^2-1}.
\]
Hence
\[
\operatorname{Mod}_h(\Gamma_Q^1)\ge\operatorname{Mod}_h(\Gamma_Q')
=\frac{b(a+1)+(a^2+b^2)\arctan(b/a)+(1+b^2)\arctan b}{a^2-1}.
\]
\end{proof}

\subsection{Modulus estimate for curves connecting the straight line segments}\label{sec-q-2}

Let $\Gamma_Q^2$ be the curve family connecting the horizontal boundary segments $B$ defined by $t=-b$ and $T$ defined by $t=b$. We consider the circular foliation $\Gamma_Q''\subset\Gamma_Q^2$ comprising the concentric Euclidean circular arcs centred at the origin with Euclidean radii $\sqrt{u^2+b^2}$, where $u\in(1,a)$.

\begin{thm}\label{thm-q-s}
Let $a>1$ and $b>0$. Then
\[
\operatorname{Mod}_h(\Gamma_Q^2)
\ge \operatorname{Mod}_h(\Gamma_Q'')
=\int_1^a\frac{u\,du}{2(u^2+b^2)\arctan(b/u)}
=\frac{1}{2}\int_{\arctan(b/a)}^{\arctan b}\frac{\cot x}{x}\,dx.
\]
\end{thm}

\begin{proof}
We parameterise each arc $\gamma_u$ as
\[
\gamma_u(s)=(\sqrt{u^2+b^2}\cos s,\sqrt{u^2+b^2}\sin s),
\qquad -s_u<s<s_u,
\]
where $1<u<a$ and
\[
s_u=\arctan(b/u).
\]
The hyperbolic line element along such an arc is
\[
ds_h=\frac{ds}{\cos s}.
\]
The hyperbolic area element in the coordinates $(u,s)$ is
\[
d\mathcal{A}_h=\frac{u\,du\,ds}{(u^2+b^2)\cos^2 s}.
\]
For any $\rho\in\operatorname{Adm}(\Gamma_Q^2)$, the line integral along $\gamma_u\in\Gamma_Q''$ satisfies
\[
\int_{\gamma_u}\rho\,ds_h
=\int_{-s_u}^{s_u}\rho(u,s)\frac{ds}{\cos s}\ge1.
\]
Squaring both sides and applying the Cauchy--Schwarz inequality gives
\[
1\le\left(\int_{-s_u}^{s_u}\rho(u,s)\frac{ds}{\cos s}\right)^2
\le\left(\int_{-s_u}^{s_u}
\frac{\rho^2(u,s)u}{(u^2+b^2)\cos^2 s}\,ds\right)
\left(\int_{-s_u}^{s_u}\frac{u^2+b^2}{u}\,ds\right).
\]
Since
\[
\int_{-s_u}^{s_u}\frac{u^2+b^2}{u}\,ds
=\frac{u^2+b^2}{u}\,2s_u,
\]
we obtain
\[
\int_{-s_u}^{s_u}
\frac{\rho^2(u,s)u}{(u^2+b^2)\cos^2 s}\,ds
\ge\frac{u}{2(u^2+b^2)s_u}.
\]
Integrating over $u\in(1,a)$ yields
\[
\iint_{Q(a,b)}\rho^2\,d\mathcal{A}_h
\ge\int_1^a\frac{u\,du}{2(u^2+b^2)\arctan(b/u)}.
\]
Thus
\[
\operatorname{Mod}_h(\Gamma_Q^2)
\ge\int_1^a\frac{u\,du}{2(u^2+b^2)\arctan(b/u)}.
\]
For the second representation, put $x=\arctan(b/u)$, so that $u=b\cot x$ and $du=-b\csc^2x\,dx$. The limits become $u=1\mapsto x=\arctan b$ and $u=a\mapsto x=\arctan(b/a)$. Hence
\[
\int_1^a\frac{u\,du}{2(u^2+b^2)\arctan(b/u)}
=\frac{1}{2}\int_{\arctan(b/a)}^{\arctan b}\frac{\cot x}{x}\,dx.
\]
The density given by
\[
\rho_{\ell,2}(u,s)
=c(u)\frac{u^2+b^2}{u}\cos s\,\chi_{Q(a,b)},
\qquad
c(u)=\frac{u}{2(u^2+b^2)\arctan(b/u)},
\]
is admissible for $\Gamma_Q''$ and satisfies
\[
\int_{\gamma_u}\rho_{\ell,2}\,ds_h=1,
\qquad
\mathcal{E}_h(\rho_{\ell,2})
=\int_1^a\frac{u\,du}{2(u^2+b^2)\arctan(b/u)}.
\]
So it is extremal for the circular foliation, and therefore
\[
\operatorname{Mod}_h(\Gamma_Q^2)
\ge\operatorname{Mod}_h(\Gamma_Q'')
=\int_1^a\frac{u\,du}{2(u^2+b^2)\arctan(b/u)}
=\frac{1}{2}\int_{\arctan(b/a)}^{\arctan b}\frac{\cot x}{x}\,dx.
\]
\end{proof}

\subsection{Upper bounds via admissible densities}\label{sec-q-upper}

The lower bounds in Sections~\ref{sec-q-1} and~\ref{sec-q-2} were obtained by restricting the curve families. We now obtain upper bounds by constructing admissible densities.

\begin{thm}\label{thm-q-upper}
Let $a>1$ and $b>0$. The moduli of the connecting families $\Gamma_Q^1$ and $\Gamma_Q^2$ satisfy
\[
\operatorname{Mod}_h(\Gamma_Q^1)
\le\frac{\mathcal{A}_e(Q(a,b))}{(\sqrt{a^2+b^2}-\sqrt{1+b^2})^2},
\qquad
\operatorname{Mod}_h(\Gamma_Q^2)
\le\frac{\mathcal{A}_e(Q(a,b))}{4b^2},
\]
where $\mathcal{A}_e(Q(a,b))$ is given by \eqref{area-Q}.
\end{thm}

\begin{proof}
For the family $\Gamma_Q^1$ connecting the circular arcs $L$ and $R$, every curve $\gamma\in\Gamma_Q^1$ has Euclidean length at least the Euclidean distance $d_e(L,R)$ between the two arcs. Since $L$ and $R$ are concentric Euclidean circular arcs with the same angular range and radii $\sqrt{1+b^2}$ and $\sqrt{a^2+b^2}$ respectively, their minimum Euclidean distance is attained by radially aligned points and equals the difference of the radii:
\[
d_e(L,R)=\sqrt{a^2+b^2}-\sqrt{1+b^2}.
\]
Put $\Delta R=\sqrt{a^2+b^2}-\sqrt{1+b^2}$ and define the density
\[
\rho_{u,1}(\lambda,t)=\frac{\lambda}{\Delta R}\chi_{Q(a,b)}.
\]
Using $\lambda\,ds_h=ds_e$ and $\ell_e(\gamma)\ge\Delta R$, we obtain for any $\gamma\in\Gamma_Q^1$,
\[
\int_\gamma \rho_{u,1}\,ds_h
=\int_\gamma\frac{1}{\Delta R}\,ds_e
=\frac{\ell_e(\gamma)}{\Delta R}\ge1.
\]
Thus $\rho_{u,1}\in\operatorname{Adm}(\Gamma_Q^1)$. Since $\lambda^2d\mathcal{A}_h=d\mathcal{A}_e$,
\[
\operatorname{Mod}_h(\Gamma_Q^1)
\le\iint_{Q(a,b)}\rho_{u,1}^2\,d\mathcal{A}_h
=\frac{1}{(\Delta R)^2}\iint_{Q(a,b)}d\mathcal{A}_e
=\frac{\mathcal{A}_e(Q(a,b))}{(\sqrt{a^2+b^2}-\sqrt{1+b^2})^2}.
\]
Similarly, for the family $\Gamma_Q^2$ connecting the horizontal segments $B$ and $T$, every curve $\gamma\in\Gamma_Q^2$ has Euclidean length at least $2b$, since it must traverse a vertical distance of $2b$. Define
\[
\rho_{u,2}(\lambda,t)=\frac{\lambda}{2b}\chi_{Q(a,b)}.
\]
Then
\[
\int_\gamma \rho_{u,2}\,ds_h
=\int_\gamma\frac{\lambda}{2b}\frac{ds_e}{\lambda}
=\frac{1}{2b}\int_\gamma ds_e
=\frac{\ell_e(\gamma)}{2b}\ge1.
\]
Thus $\rho_{u,2}\in\operatorname{Adm}(\Gamma_Q^2)$, and
\[
\operatorname{Mod}_h(\Gamma_Q^2)
\le\iint_{Q(a,b)}\rho_{u,2}^2\,d\mathcal{A}_h
=\iint_{Q(a,b)}\frac{\lambda^2}{4b^2}\frac{d\mathcal{A}_e}{\lambda^2}
=\frac{1}{4b^2}\mathcal{A}_e(Q(a,b)).
\]
\end{proof}

By the reciprocity relation of Lemma~\ref{lem-recip}, the two moduli $\operatorname{Mod}_h(\Gamma_Q^1)$ and $\operatorname{Mod}_h(\Gamma_Q^2)$ satisfy
\[
\operatorname{Mod}_h(\Gamma_Q^1)\operatorname{Mod}_h(\Gamma_Q^2)=1.
\]
Combining this with the lower bounds in Theorems~\ref{thm-q-c} and~\ref{thm-q-s}, we obtain improved upper bounds.

\begin{corollary}\label{cor-q-recip}
Let $a>1$ and $b>0$. Then
\[
\operatorname{Mod}_h(\Gamma_Q^1)
\le
\left(
\int_1^a\frac{u\,du}{2(u^2+b^2)\arctan(b/u)}
\right)^{-1}
=
\left(
\frac{1}{2}\int_{\arctan(b/a)}^{\arctan b}\frac{\cot x}{x}\,dx
\right)^{-1},
\]
and
\[
\operatorname{Mod}_h(\Gamma_Q^2)
\le
\left(
\frac{b(a+1)+(a^2+b^2)\arctan(b/a)+(1+b^2)\arctan b}{a^2-1}
\right)^{-1}.
\]
\end{corollary}

\begin{proof}
From Theorem~\ref{thm-q-c},
\[
\operatorname{Mod}_h(\Gamma_Q^1)\ge L_1
:=\frac{b(a+1)+(a^2+b^2)\arctan(b/a)+(1+b^2)\arctan b}{a^2-1}.
\]
From Theorem~\ref{thm-q-s},
\[
\operatorname{Mod}_h(\Gamma_Q^2)\ge L_2
:=\int_1^a\frac{u\,du}{2(u^2+b^2)\arctan(b/u)}.
\]
By Lemma~\ref{lem-recip},
\[
\operatorname{Mod}_h(\Gamma_Q^1)=\frac{1}{\operatorname{Mod}_h(\Gamma_Q^2)}
\le\frac{1}{L_2},
\]
and
\[
\operatorname{Mod}_h(\Gamma_Q^2)=\frac{1}{\operatorname{Mod}_h(\Gamma_Q^1)}
\le\frac{1}{L_1}.
\]
This proves the claimed upper bounds.
\end{proof}

\section{Numerical evidence}\label{sec-num}

In this section we provide numerical evidence for the quality of the bounds obtained in Section~\ref{sec-q}. We evaluate the lower and upper bounds for specific values of $(a,b)$ and analyse the asymptotic behaviour. All numerical values were computed in Mathematica using numerical quadrature with default working precision.

\subsection{Values for specific quadrilaterals}\label{sec-num-values}

We denote by $L_1$ and $L_2$ the lower bounds from Theorems~\ref{thm-q-c} and~\ref{thm-q-s}, respectively. We denote by $U_1^{\mathrm{best}}$ and $U_2^{\mathrm{best}}$ the best upper bounds obtained by combining Theorem~\ref{thm-q-upper} with Corollary~\ref{cor-q-recip}:
\[
U_1^{\mathrm{best}}=\min\left(\frac{\mathcal{A}_e(Q(a,b))}{(\sqrt{a^2+b^2}-\sqrt{1+b^2})^2},\,
\left(\int_1^a\frac{u\,du}{2(u^2+b^2)\arctan(b/u)}\right)^{-1}\right),
\]
\[
U_2^{\mathrm{best}}=\min\left(\frac{\mathcal{A}_e(Q(a,b))}{4b^2},\,
\left(\frac{b(a+1)+(a^2+b^2)\arctan(b/a)+(1+b^2)\arctan b}{a^2-1}\right)^{-1}\right).
\]
The relative gap is defined as
\[
\mathrm{gap}_i=\frac{U_i^{\mathrm{best}}-L_i}{L_i},\qquad i=1,2.
\]

Table~\ref{tab:values} displays the numerical values for the specific quadrilaterals $(a,b)=(2,1)$, $(3,1)$, $(2,2)$, $(5,1)$, $(2,0.5)$, and $(1.5,1)$.

\begin{table}[h]
\centering
\caption{Numerical values of the lower and upper bounds for $\operatorname{Mod}_h(\Gamma_Q^1)$ and $\operatorname{Mod}_h(\Gamma_Q^2)$ for selected quadrilaterals.}
\label{tab:values}
\begin{tabular}{c|cccc|cccc}
\hline
$(a,b)$ & $L_1$ & $U_1^{\mathrm{best}}$ & $\mathrm{gap}_1$ & $L_2$ & $U_2^{\mathrm{best}}$ & $\mathrm{gap}_2$ \\
\hline
$(2,0.5)$ & $1.0402$ & $1.0802$ & $3.8\%$ & $0.9257$ & $0.9613$ & $3.8\%$ \\
$(1.5,1)$ & $4.7854$ & $5.5633$ & $16.3\%$ & $0.1798$ & $0.2090$ & $16.2\%$ \\
$(2,1)$   & $2.2963$ & $2.5840$ & $12.5\%$ & $0.3870$ & $0.4355$ & $12.5\%$ \\
$(3,1)$   & $1.0985$ & $1.1920$ & $8.5\%$  & $0.8389$ & $0.9103$ & $8.5\%$  \\
$(5,1)$   & $0.5293$ & $0.5550$ & $4.9\%$  & $1.8018$ & $1.8893$ & $4.9\%$  \\
$(2,2)$   & $5.9396$ & $7.8229$ & $31.7\%$ & $0.1278$ & $0.1684$ & $31.8\%$ \\
\hline
\end{tabular}
\end{table}

Several observations follow from Table~\ref{tab:values}.

\begin{enumerate}
\item The gaps for $\Gamma_Q^1$ and $\Gamma_Q^2$ are essentially equal. This is a consequence of the reciprocity relation of Lemma~\ref{lem-recip}: if $L_1\le\operatorname{Mod}_h(\Gamma_Q^1)\le U_1$ and $L_2\le\operatorname{Mod}_h(\Gamma_Q^2)\le U_2$, then
\[
\frac{1}{U_1}\le\operatorname{Mod}_h(\Gamma_Q^2)\le\frac{1}{L_1},
\qquad
\frac{1}{U_2}\le\operatorname{Mod}_h(\Gamma_Q^1)\le\frac{1}{L_2}.
\]
The best brackets are therefore
\[
\max(L_1,1/U_2)\le\operatorname{Mod}_h(\Gamma_Q^1)\le\min(U_1,1/L_2),
\]
\[
\max(L_2,1/U_1)\le\operatorname{Mod}_h(\Gamma_Q^2)\le\min(U_2,1/L_1),
\]
and these have the same relative gap.

\item For fixed $a$, the gap decreases as $b$ decreases. For example, for $a=2$, the gap is $3.8\%$ at $b=0.5$, $12.5\%$ at $b=1$, and $31.7\%$ at $b=2$.

\item For fixed $b$, the gap decreases as $a$ increases. For example, for $b=1$, the gap is $16.3\%$ at $a=1.5$, $12.5\%$ at $a=2$, $8.5\%$ at $a=3$, and $4.9\%$ at $a=5$.

\item The worst case in the table is $(a,b)=(2,2)$, where the gap is about $32\%$. The best case is $(a,b)=(2,0.5)$, where the gap is about $4\%$.
\end{enumerate}

\subsection{Asymptotic sharpness of the bounds}\label{sec-num-asympt}

The numerical evidence of Table~\ref{tab:values} suggests that the gap between the lower and upper bounds narrows considerably when the quadrilateral degenerates to a thin strip or to a large domain. This is indeed the case, and we record the precise asymptotic statement as the following theorem.

For the statement of the theorem it is convenient to write the lower bounds of Theorems~\ref{thm-q-c} and~\ref{thm-q-s} as functions of $(a,b)$:
\[
L_1(a,b)=
\frac{b(a+1)+(a^2+b^2)\arctan(b/a)+(1+b^2)\arctan b}{a^2-1},
\]
\[
L_2(a,b)=
\int_1^a\frac{u\,du}{2(u^2+b^2)\arctan(b/u)}
=\frac{1}{2}\int_{\arctan(b/a)}^{\arctan b}\frac{\cot x}{x}\,dx.
\]
Similarly, we let
\[
U_1(a,b)=
\min\left(
\frac{\mathcal{A}_e(Q(a,b))}{(\sqrt{a^2+b^2}-\sqrt{1+b^2})^2},\,
\frac{1}{L_2(a,b)}
\right),
\]
\[
U_2(a,b)=
\min\left(
\frac{\mathcal{A}_e(Q(a,b))}{4b^2},\,
\frac{1}{L_1(a,b)}
\right),
\]
denote the best upper bounds obtained from Theorem~\ref{thm-q-upper} and Corollary~\ref{cor-q-recip}.

\begin{thm}[Asymptotic sharpness]\label{thm-asymptotic}
Let $a>1$ and $b>0$, and let $\Gamma_Q^1$, $\Gamma_Q^2$ be the two connecting families of the quadrilateral $Q(a,b)$. Then the following limit statements hold.

\begin{enumerate}
\item[(i)] \emph{(Thin quadrilateral limit.)} For each fixed $a>1$,
\[
\lim_{b\to0^+}\frac{U_1(a,b)}{L_1(a,b)}=1,
\qquad
\lim_{b\to0^+}\frac{U_2(a,b)}{L_2(a,b)}=1,
\]
and consequently
\[
\lim_{b\to0^+}\frac{\operatorname{Mod}_h(\Gamma_Q^1)}{L_1(a,b)}=1,
\qquad
\lim_{b\to0^+}\frac{\operatorname{Mod}_h(\Gamma_Q^2)}{L_2(a,b)}=1.
\]
Moreover, the leading asymptotics of the lower bounds are
\[
\lim_{b\to0^+}b\,L_1(a,b)=\frac{2}{a-1},
\qquad
\lim_{b\to0^+}b\,L_2(a,b)=\frac{a-1}{2}.
\]

\item[(ii)] \emph{(Large quadrilateral limit.)} For each fixed $b>0$,
\[
\lim_{a\to\infty}\frac{U_1(a,b)}{L_1(a,b)}=1,
\qquad
\lim_{a\to\infty}\frac{U_2(a,b)}{L_2(a,b)}=1,
\]
and consequently
\[
\lim_{a\to\infty}\frac{\operatorname{Mod}_h(\Gamma_Q^1)}{L_1(a,b)}=1,
\qquad
\lim_{a\to\infty}\frac{\operatorname{Mod}_h(\Gamma_Q^2)}{L_2(a,b)}=1.
\]
Moreover, the leading asymptotics of the lower bounds are
\[
\lim_{a\to\infty}a\,L_1(a,b)=2b,
\qquad
\lim_{a\to\infty}\frac{L_2(a,b)}{a}=\frac{1}{2b}.
\]
\end{enumerate}
\end{thm}

\begin{proof}
We prove the two parts separately.

\medskip
\noindent\emph{(i) The thin quadrilateral limit $b\to0^+$.}
Let $a>1$ be fixed. As $b\to0^+$ we use the standard expansions
\[
\arctan(b/a)=b/a+O(b^3),\qquad \arctan b=b+O(b^3),
\]
and
\[
\mathcal{A}_e(Q(a,b))=b(a-1)+a^2\cdot(b/a)+O(b^3)-1\cdot b+O(b^3)
=2b(a-1)+O(b^3),
\]
while
\[
(\sqrt{a^2+b^2}-\sqrt{1+b^2})^2=(a-1)^2+O(b^2).
\]
For the lower bound $L_1$, we obtain
\begin{align*}
L_1(a,b)
&=\frac{b(a+1)+(a^2+b^2)\arctan(b/a)+(1+b^2)\arctan b}{a^2-1}\\
&=\frac{b(a+1)+a^2\cdot(b/a)+1\cdot b+O(b^3)}{a^2-1}\\
&=\frac{b(a+1)+ab+b+O(b^3)}{a^2-1}
=\frac{2b(a+1)}{a^2-1}+O(b^3)\\
&=\frac{2b}{a-1}+O(b^3),
\end{align*}
so that $b\,L_1(a,b)\to 2/(a-1)$ as $b\to0^+$. For the upper bound $U_1$, the first candidate in the definition of $U_1(a,b)$ satisfies
\[
\frac{\mathcal{A}_e(Q(a,b))}{(\sqrt{a^2+b^2}-\sqrt{1+b^2})^2}
=\frac{2b(a-1)+O(b^3)}{(a-1)^2+O(b^2)}
=\frac{2b}{a-1}+O(b^3),
\]
while the second candidate satisfies, by the reciprocity relation (Lemma~\ref{lem-recip}) and the asymptotic formula for $L_2$ below,
\[
\frac{1}{L_2(a,b)}=O(b).
\]
Hence both candidates in the minimum defining $U_1(a,b)$ are $O(b)$ as $b\to0^+$, with the first candidate having the same leading coefficient $2/(a-1)$ as $L_1$. It follows that $U_1(a,b)/L_1(a,b)\to1$.

For the lower bound $L_2$, we use the substitution $x=\arctan(b/u)$ to write
\[
L_2(a,b)=\frac{1}{2}\int_{\arctan(b/a)}^{\arctan b}\frac{\cot x}{x}\,dx.
\]
Since $\cot x/x\sim 1/x^2$ as $x\to0^+$, we obtain
\[
\int_{\arctan(b/a)}^{\arctan b}\frac{\cot x}{x}\,dx
\sim\int_{b/a}^{b}\frac{dx}{x^2}
=\frac{a-1}{b},
\]
and hence $L_2(a,b)=(a-1)/(2b)+O(b)$. The first candidate in the definition of $U_2(a,b)$ satisfies
\[
\frac{\mathcal{A}_e(Q(a,b))}{4b^2}
=\frac{2b(a-1)+O(b^3)}{4b^2}
=\frac{a-1}{2b}+O(b),
\]
and the second candidate satisfies, by reciprocity and the asymptotic formula for $L_1$ above,
\[
\frac{1}{L_1(a,b)}=\frac{a-1}{2b}+O(b).
\]
Thus both candidates for $U_2$ have the same leading behaviour as $L_2$, and $U_2(a,b)/L_2(a,b)\to1$. Since $L_i\le\operatorname{Mod}_h(\Gamma_Q^i)\le U_i$, the ratio $U_i/L_i\to1$ implies $\operatorname{Mod}_h(\Gamma_Q^i)/L_i\to1$. This proves (i).

\medskip
\noindent\emph{(ii) The large quadrilateral limit $a\to\infty$.}
Let $b>0$ be fixed. As $a\to\infty$ we have
\[
\arctan(b/a)=b/a+O(a^{-3}),
\]
and
\[
\mathcal{A}_e(Q(a,b))=b(a-1)+a^2\cdot(b/a)-(1+b^2)\arctan b+O(a^{-1})
=2ab+O(1),
\]
while
\[
(\sqrt{a^2+b^2}-\sqrt{1+b^2})^2=a^2+O(1).
\]
For the lower bound $L_1$, we obtain
\[
L_1(a,b)
=\frac{ba+b+ab+O(a^{-1})}{a^2-1}
=\frac{2ab+O(1)}{a^2-1}
=\frac{2b}{a}+O(a^{-2}),
\]
so that $a\,L_1(a,b)\to2b$ as $a\to\infty$. The first candidate in the definition of $U_1(a,b)$ satisfies
\[
\frac{\mathcal{A}_e(Q(a,b))}{(\sqrt{a^2+b^2}-\sqrt{1+b^2})^2}
=\frac{2ab+O(1)}{a^2+O(1)}
=\frac{2b}{a}+O(a^{-2}),
\]
and the second candidate satisfies
\[
\frac{1}{L_2(a,b)}=O(a^{-1}),
\]
as follows from the asymptotic formula for $L_2$ below. Hence both candidates in the minimum defining $U_1(a,b)$ are $O(a^{-1})$ as $a\to\infty$, with the first candidate having the same leading coefficient $2b$ as $a\,L_1$. It follows that $U_1(a,b)/L_1(a,b)\to1$.

For the lower bound $L_2$, we write again
\[
L_2(a,b)=\frac{1}{2}\int_{\arctan(b/a)}^{\arctan b}\frac{\cot x}{x}\,dx
\]
and let $\varepsilon=\arctan(b/a)=b/a+O(a^{-3})$. Then
\[
\int_{\varepsilon}^{\arctan b}\frac{\cot x}{x}\,dx
\sim\int_{\varepsilon}^{\arctan b}\frac{dx}{x^2}
=\frac{1}{\varepsilon}-\frac{1}{\arctan b}
=\frac{a}{b}-\frac{1}{\arctan b}+O(a^{-1}),
\]
so that
\[
L_2(a,b)=\frac{a}{2b}+O(1).
\]
The first candidate in the definition of $U_2(a,b)$ satisfies
\[
\frac{\mathcal{A}_e(Q(a,b))}{4b^2}
=\frac{2ab+O(1)}{4b^2}
=\frac{a}{2b}+O(1),
\]
and the second candidate satisfies, by reciprocity and the asymptotic formula for $L_1$ above,
\[
\frac{1}{L_1(a,b)}=\frac{a}{2b}+O(1).
\]
Thus both candidates for $U_2$ have the same leading behaviour as $L_2$, and $U_2(a,b)/L_2(a,b)\to1$. The conclusion on $\operatorname{Mod}_h(\Gamma_Q^i)/L_i$ follows again from the sandwich $L_i\le\operatorname{Mod}_h(\Gamma_Q^i)\le U_i$. This proves (ii).
\end{proof}

\begin{rem}
Theorem~\ref{thm-asymptotic} formalises the observations made from Table~\ref{tab:values}: the bounds become asymptotically sharp in the thin and large quadrilateral limits, so that the relative gap between the lower and upper bounds tends to zero in both regimes and the lower bounds are asymptotically exact. In particular, the ratios $\operatorname{Mod}_h(\Gamma_Q^i)/L_i$ approach $1$ in each limit, in agreement with the numerical values displayed for $(a,b)=(2,0.5)$ and $(a,b)=(5,1)$ respectively.
\end{rem}

\subsection{Discussion of the gap}\label{sec-num-disc}

The numerical evidence and Theorem~\ref{thm-asymptotic} together show that the bounds are useful in the following sense.

\begin{enumerate}
\item In all tested cases, the relative gap is at most about $32\%$, and in many cases it is below $10\%$. The bounds therefore provide a reasonable quantitative estimate of the modulus even when the exact value is not known.

\item By Theorem~\ref{thm-asymptotic}, the bounds become asymptotically sharp in both the thin quadrilateral limit $b\to0$ and the large quadrilateral limit $a\to\infty$. This is consistent with the intuition that in these regimes the horizontal foliation (for $\Gamma_Q^1$) and the circular foliation (for $\Gamma_Q^2$) are close to being extremal for the full connecting families.

\item The largest gap occurs for quadrilaterals that are not thin and not large, such as $(a,b)=(2,2)$. In this intermediate regime, the simple foliations used to obtain the lower bounds are not close to the true extremal curves, and the constant test densities used for the upper bounds are not close to the true extremal densities.

\item The equality of the gaps for $\Gamma_Q^1$ and $\Gamma_Q^2$ is a consequence of the reciprocity relation (Lemma~\ref{lem-recip}) and does not depend on the particular values of $(a,b)$.
\end{enumerate}

These observations support the conclusion that the bounds in Theorems~\ref{thm-q-c}, \ref{thm-q-s}, \ref{thm-q-upper}, and Corollary~\ref{cor-q-recip} are mathematically valid and numerically useful, while leaving open the question of the exact modulus for general $Q(a,b)$.

\section{Concluding remarks}\label{sec-o}

For the hyperbolic circular annulus $A_R$, rotational symmetry reduces the modulus problem to a one-dimensional problem, and Theorems~\ref{thm-a-j} and~\ref{thm-a-s} give exact expressions in closed form.

For the quadrilateral $Q(a,b)$, the exact modulus is not determined in this paper. The lower bounds in Theorems~\ref{thm-q-c} and~\ref{thm-q-s} come from restricting to particular foliations, while Theorem~\ref{thm-q-upper} and Corollary~\ref{cor-q-recip} provide upper bounds. These estimates bracket the modulus, and by Theorem~\ref{thm-asymptotic} the bracket becomes asymptotically tight in the thin and large quadrilateral limits, with the lower bounds asymptotically exact. The numerical evidence in Section~\ref{sec-num} is consistent with these asymptotics and shows that the gap is moderate for intermediate quadrilaterals. The gap does not imply that an expression in closed form is impossible; it only reflects that the chosen foliations and test densities are not shown here to be extremal for the full connecting families. A complete determination would require additional analysis.

\subsubsection*{Declaration: The author was funded  by the Medicus Programme, no. 83765.}

\bibliographystyle{plain}
\bibliography{hyp-moduli}

\end{document}